\documentstyle[12pt]{article} 
\newfont{\twlvmsb}{msbm10 scaled\magstep1}
\newfont{\ninemsb}{msbm9}
\newfont{\sixmsb}{msbm6}
\newfam\msbfam
\textfont\msbfam=\twlvmsb
\scriptfont\msbfam=\ninemsb
\scriptscriptfont\msbfam=\sixmsb
\catcode`\@=11
\def\Bbb{\ifmmode\let\next\Bbb@\else
  \def\next{\errmessage{Use \string\Bbb\space only in math mode}}\fi\next}
\def\Bbb@#1{{\Bbb@@{#1}}}
\def\Bbb@@#1{\fam\msbfam#1}
\newfont{\largeeufm}{eufm10 scaled\magstep4}
\newfont{\twlveufm}{eufm10 scaled\magstep1}
\newfont{\elveufm}{eufm10 at 11pt}
\newfont{\teneufm}{eufm10}
\newfont{\nineeufm}{eufm9}
\newfam\eufam
\textfont\eufam=\twlveufm
\scriptfont\eufam=\nineeufm
\def\frak{\ifmmode\let\next\frak@\else
\def\next{\errmessage{Use \string\frak\space only in math mode}}\fi\next}
\def\frak@#1{{\fam\eufam{{#1}}}}
\catcode`@=12
\newcommand{\C}{{\Bbb C}} 
\newcommand{\g}{{\frak g}} 
\newcommand{\Uq}{U_q({\frak g})} 
\newcommand{\U}{U({\frak g})} 
\newcommand{\Up}{U_q^{+}} 
\newcommand{\Um}{U_q^{-}} 
\newcommand{\ba}{\begin{eqnarray}}
\newcommand{\na}{\end{eqnarray}}
\newcommand{\ban}{\begin{eqnarray*}}
\newcommand{\nan}{\end{eqnarray*}}

\begin{document} 
\title{\normalsize{\bf MINIMAL UNCERTAINTY STATES FOR QUANTUM GROUPS} } 
\author{ \small R. DELBOURGO$^1$ and R. B. ZHANG$^2$ \\  
\small $^1$Department of Physics, University of Tasmania,
Hobart, Australia\\ 
\small $^2$Department of Pure Mathematics,  University of Adelaide, 
Adelaide, Australia}

\date{} 
\maketitle
\begin{abstract} 
The problem of how to
obtain quasi-classical states for quantum groups is examined. 
A measure of quantum indeterminacy is proposed, which 
involves expectation values of some natural quantum group operators. 
It is shown that within any finite dimensional irreducible representation, 
the highest weight vector and those unitarily related to it 
are the quasi-classical states.  
\end{abstract}

Quantum groups \cite{Drinfeld} have been intensively studied in  
recent years. Their applications  
have already led to significant progress in statistical 
mechanics and low dimensional topology. It is also widely 
believed that quantum groups play some important role in 
quantum physics as well.  
In this letter, we will investigate the problem of how to
obtain quasi classical states for quantum groups. 
We will propose a measure of quantum indeterminacy, 
which involves expectation values of some combinations of 
Drinfeld's $v$ operator and the universal $R$ - matrix. 
A quasi-classical state is characterized as having 
minimal indeterminacy. It will be shown that 
for any finite dimensional irreducible representation, 
the highest weight vector and those unitarily related to 
it are the states having this property.  
Our study here is an extension to quantum 
groups of the investigation carried out in \cite{Delbourgo} 
some twenty years ago, 
where the corresponding problem for compact simple Lie groups 
was resolved by one of us. 
In the limit $q\rightarrow 1$, we recover the results of 
that publication.

Given a compact simple Lie group $G$, we denote its 
Lie algebra by $\g$.  Now there exists a
basis $\{ e_i\}$, which is self dual with respect to the Killing 
form, in which the quadratic Casimir operator can  
be expressed as $C=\sum_i e_i \, e_i$.   
It was shown in \cite{Delbourgo} that the following dispersion 
\ba
\sum_i \langle \left( e_i - \langle e_i \rangle\right) 
\left( e_i - \langle e_i \rangle\right) \rangle \label{dispersion}  
\na 
was the natural measure of quantum indeterminacy for such an
algebra. For any finite dimensional irreducible representation, 
the highest weight vector and its group orbit, i.e., 
the Peremolov coherent states \cite{Peremolov}, 
proved to be the vectors attaining minimal uncertainty.  

Let us rewrite (\ref{dispersion}) in a more abstract form, 
so that it will suggests a  generalization to quantum groups.     
Acting on any irreducible representation space with 
a highest weight $\lambda$, $C$ takes the eigenvalue 
$(\lambda + 2\rho,\ \lambda)$, where $2\rho$ is the sum of all 
the positive roots, and $( .\ ,  . )$ is the properly 
normalized inner product of the weight space.  However, when 
acting on the tensor product of two representation spaces, 
$C$ will no longer be a scalar multiple of the identity matrix. 
We  express the action of $C$ on the tensor product by 
$\partial (C)$, and set  
$Q=\left[ \partial (C) - C\otimes 1 - 1\otimes C\right]$.     
For any unit vector $|\psi\rangle$, the indeterminacy 
measure (\ref{dispersion}) can be expressed as   
\ba 
\langle\psi| C|\psi \rangle - 
\langle \psi\otimes\psi| Q |\psi\otimes\psi \rangle,\label{measure}  
\na 
where we have used the self-explanatory notation 
that $|\psi\otimes\psi \rangle = |\psi\rangle \otimes|\psi\rangle$.

Recall that a quantum group $\Uq$ is a deformation of the universal 
enveloping algebra $\U$ of the Lie algebra $\g$. It has the 
structures of a Hopf algebra, namely, it admits a tensor product 
operation, called the co-multiplication, 
 $\Delta: \Uq\rightarrow \Uq\otimes\Uq$, a co-unit $\epsilon: \Uq
\rightarrow\C$, and an antipode $S: \Uq\rightarrow \Uq$, which 
are all compatible with the algebraic structure of $\Uq$ in an 
appropriate sense.  $\Uq$ also has properties very similar to those of  
$\U$, e.g., it is  generated by raising, lowering and 
Cartan type of generators.  The deformation parameter may be 
assumed to take various forms, leading to different versions of 
quantum groups.  Here we will take $q=\exp(\eta)$ with 
$\eta$ being a real number not equal to $1$. 
For the purpose of this letter, we will allow exponentials of 
Cartan generators and their polynomials to appear in the algebra 
as well. Since we will limit ourselves to finite dimensional 
representations, such things make perfect sense. 

We denote by $\Up$ the subalgebra generated by the raising and Cartan 
elements, and by $\Um$ that generated  by the lowering and Cartan elements. 
There exist bases $\{\alpha_t\, | \, t=1, 2, ... \}$ 
and $\{\beta_t\, | \, t=1, 2, ... \}$ for $\Up$ and $\Um$ 
respectively such that 
\ban 
R&=& \sum_t \alpha_t\otimes \beta_t 
\nan 
gives rise to the universal $R$ - matrix of $\Uq$.  
Furthermore, $\Uq$ admits an involution $ ^\dagger$ satisfying 
$ ^\dagger\, S\,  ^\dagger\, S= 1$, rendering 
$\alpha_t^\dagger =\pm  \beta_t$.   With respect to this involution 
all finite dimensional representations are unitary.

Let $K_\rho$ be the group like element of $\Uq$ such 
that  $S^2(x)=K_\rho^2 \, x\, K_\rho^{-2}$ for all $x$ in $\Uq$.      
The Drinfeld operator of the quantum group is given by 
$v= \sum_t S(\beta_t)\, \alpha_t\ K_\rho^{-2}$,  
which is central, namely, commutes with all the
elements of $\Uq$. It is also known that 
$v$ is fixed by the antipode, i.e., $S(v)=v$, 
and is invertible with the inverse given by  
$v^{-1} =\sum_t \beta_t \, K_\rho^2 \, \alpha_t$. 
In a finite dimensional irreducible representation $V(\lambda)$ 
with highest weight $\lambda$,  the operator $v^{-1}$ is 
given by $q^{(\lambda + 2\rho,\ \lambda)}\ I$.    

In studying the quasi-classical states, we will need the operator 
$v^{-2}$, which can be expressed as   
\ban
v^{-2} &=& \sum_{r,\, t} \beta_r\, S(\alpha_t)\,  
S^{-1}(\beta_t)\,  \alpha_r.
\nan     
We will also need the operator  
 $(v\otimes v) \Delta(v^{-1}) = R^T R$,  where  
\ba 
R^T R &=& \sum_{r,\, t} \beta_r\, \alpha_t\, \otimes  
\alpha_r\, \beta_t\nonumber\\  
&=& \sum_{r,\, t} \beta_r\, S(\alpha_t)\, \otimes \alpha_r\, 
S(\beta_t). \label{tensor} 
\na 

$R^T R$ acts naturally on the tensor product of two 
representation spaces. Consider for example the tensor product 
of $V(\lambda)$ with 
itself, which can always be decomposed into a direct sum of finite 
dimensional irreducible representations   
\ban V(\lambda)\otimes V(\lambda) = \bigoplus_{i=0}^L V(\mu_i), \nan 
where $L$ is some positive integer  which we will 
not need to know.     The $V(\mu_i)$ is an irreducible representation 
with highest weight $\mu_i$. We will order the $\mu$'s in such a way 
that $\mu_i\ge \mu_j$ if $i<j$. Then clearly, $\mu_0=2\lambda$ $>\mu_i$, 
for all $i>0$. 
$R^T R$ in $V(\lambda)\otimes V(\lambda)$, though not being 
proportional to the identity matrix,  can nevertheless 
be diagonalized, and its eigenvalues are 
$$q^{(\mu_i+2\rho,\, \mu_i) - 
2(\lambda + 2\rho,\ \lambda)},\qquad i=0, 1, ..., L.$$  
 
Let $|\phi\rangle$ be a unit vector of $V(\lambda)$. 
The following quantity is a natural generalization of (\ref{measure}) 
to the quantum group setting: 
\ba 
\delta_\phi &=& { {1}\over{q-q^{-1}} } 
\left[ \langle\phi| v^{-2}|\phi\rangle  
- \langle\phi\otimes\phi| R^T R|\phi\otimes\phi\rangle\right].  
\na 
Observe that both terms of $\delta_\phi$  are invariant with 
respect to the quantum group, i.e., 
\ban 
 \langle\phi| [x, \  v^{-2}] |\phi\rangle&=&0, \\ 
\langle\phi\otimes\phi| [\Delta(x), \  R^T R] |\phi\otimes\phi\rangle
&=&0,   \ \ \ \forall x\in \Uq.  
\nan   
Also,  we recover from $\delta_\phi$ the dispersion (\ref{dispersion}) 
for Lie algebras in the limit $q\rightarrow 1$.  
We propose  $\delta_\phi$  
as the measure of quantum indeterminacy, and will call 
a quantum group state quasi-classical if $\delta_\phi$ is minimized. 
 
In order to minimize $\delta_\phi$, we observe that the 
vector $|\phi\otimes\phi\rangle$ can always be written  as 
\ban 
|\phi\otimes\phi\rangle&=&\sum_{i=0}^L c_i|\zeta_i\rangle\rangle,
\nan  
where $|\zeta_i\rangle\rangle$ is a unit vector belonging to 
$V(\mu_i)$, and   $c_i$ is a {\em real number} satisfying
the normalization property 
$\sum_{i=0}^L c_i^2=1$.  With the help of this expression, 
we can  easily obtain   
\ban \delta_\phi &=& { {1}\over{q-q^{-1}} }
\left[q^{2(\lambda + 2\rho,\ \lambda)} 
- \sum_{i=0}^L c_i^2\, q^{(\mu_i+2\rho,\ \mu_i) -
2(\lambda + 2\rho,\ \lambda)}\right]. 
\nan 
Since $\mu_i -\mu_j$, $i<j$,  is either zero or 
a positive integral sum of positive roots of $\g$, we have  
\ban 
(\mu_i+2\rho,\ \mu_i)- (\mu_j+2\rho,\ \mu_j) 
&=&(\mu_i+\mu_j+2\rho,\ \mu_i-\mu_j)\ge 0, \ \ \ i<j.    
\nan 
This immediately leads to the  conclusion that   
$\delta_\phi$ reaches its minimum  
\ban 
Min(\delta)&=&{ {1}\over{q-q^{-1}} } 
\left[q^{2(\lambda + 2\rho,\ \lambda)}
- q^{2(\lambda,\ \lambda)}\right]\ge 0,    
\nan   
only when 
\ban c_0^2=1,  &c_1=...=c_L=0. \nan   
Therefore, {\em  a state  $|\phi\rangle\in V(\lambda)$
is quasi-classical if and only if $|\phi\otimes\phi\rangle$ 
belongs to the irreducible component $V(2\lambda)$ contained in 
$V(\lambda)\otimes V(\lambda)$}. 
It is clearly true that the highest weight vector $|\phi_0\rangle$ 
(normalized to $1$) of $V(\lambda)$ meets this requirement.  
Given any unit vector $|\phi\rangle$,  the irreducibility of 
$V(\lambda)$ guarantees the existence of a unitary endomorphism   
$U$ ( i.e., $U^\dagger \, U = 1$ ) such that 
$|\phi\rangle=U|\phi_0\rangle$. When $U\otimes U$ 
commutes with $R^T R$, $|\phi\rangle$ yields a  
quasi-classical state, and all the quasi-classical states 
are of this form. In the limit $q\rightarrow 1$, these states 
reduce to Peremolov's  coherent states. 

A deeper understanding of the quasi-classical states, 
particularly their underlying  geometry,   could be 
gained by studying their properties with respect to the 
algebra of functions on $\Uq$, but this would take us into 
the largely unexplored area of noncommutative geometry, 
which is well beyond the scope of this letter.

\vspace{1cm}
\noindent
{\bf Acknowledgement}: We thank Drs. K. Hannabuss and P. Jarvis 
for helpful discussions.

 
\end{document}